\documentclass{article}

\usepackage[all]{xy}
 \usepackage[usenames,dvipsnames]{pstricks}
 \usepackage{epsfig}
 \usepackage{pst-grad} % For gradients
 \usepackage{pst-plot} % For axes
\usepackage{ebezier}

\usepackage[
        colorlinks, citecolor=red,
        backref,
        pdfauthor={MDJ,JK},
        pdftitle={S1}
]{hyperref}

\usepackage{latexsym}\usepackage{amsfonts}\usepackage{amssymb}\usepackage{amsthm}
\usepackage{amsmath}\usepackage{newlfont}\usepackage{graphicx}\usepackage{doc}
\usepackage{color}\usepackage{multicol}\usepackage{epic}\usepackage{eepic}
\usepackage{ebezier}
\newtheorem{thrm}{Theorem}
\newtheorem{lem}[thrm]{Lemma}
\newtheorem{cor}[thrm]{Corollary}

\newtheorem{prop}[thrm]{Proposition}
\newtheorem*{theorem-non}{Theorem}

\newtheorem*{namedthm}{\namedthmname}
\newcounter{namedthm}



\include{epsf}

\def \Dj{\mbox{\raise0.3ex\hbox{-}\kern-0.4em D}}

\begin{document}

\title{Polynomial entropy of induced homeomorphisms on $C_k(\mathbb{S}^1)$ and $C_k([0,1])$}

\author{Ma\v{s}a \Dj ori\'c\thanks{Corresponding author: Ma\v{s}a \Dj ori\'c.}\\
Matemati\v{c}ki institut SANU\\
Knez Mihailova 36\\
11000 Beograd\\
Serbia\\
masha@mi.sanu.ac.rs\\ \and
Jelena Kati\'c\\
Matemati\v cki fakultet\\
Studentski trg 16\\
11000 Beograd\\
Serbia\\
jelena.katic@matf.bg.ac.rs
}

\maketitle

\begin{abstract} We compute the polynomial entropy of $C_k(f)$ where $f$ is any circle or interval homeomorphism.
\end{abstract}

\medskip

\section{Introduction}

Every continuous map on a compact metric space $X$ induces a continuous map (called the \textit{induced map}) on the hyperspace $2^X$ of all closed subsets. If $X$ is connected, i.e.\ a continuum, we consider the hyperspace $C(X)$ of subcontinua of $X$ which is itself a continuum. The space $C_k(X)$ of all nonempty compact subsets of $X$ having at most $k$ connected components is another well-studied hyperspace. A natural question is: what are the possible relations between the given (individual) dynamics on $X$ and the induced one (collective) dynamics on the hyperspace. In recent decades, various results have been obtained in this direction. Without claiming to be exhaustive, we will mention just a few: Borsuk and Ulam~\cite{BU}, Bauer and Sigmund~\cite{BS}, Rom\'an-Flores~\cite{RF}, Banks~\cite{B}, Acosta, Illanes and M\'endez-Lango~\cite{AIM}.

The topological entropy of the induced map has been studied by Kwietniak and Oprocha in~\cite{KO}, Lampart and Raith~\cite{LR}, Hern\'andez and M\'endez~\cite{HM}, Arbieto and Bohorquez~\cite{AB}, and others.

In~\cite{LR} the authors showed that, if $f$ is an interval or a circle homeomorphism, the topological entropy of the induced map on the hyperspace of subcontinua is zero. This result also appears as a corollary in~\cite{AB}, in the context of Morse–Smale diffeomorphisms of the circle.

One useful measure of complexity for systems with zero topological entropy is polynomial entropy. The very definitions of topological and polynomial entropy differ in that the former quantifies the overall {\it exponential} growth of orbit complexity, while the latter captures the growth rate on a {\it polynomial} scale. Polynomial entropy serves as a refined tool to distinguish among systems with zero topological entropy. For instance, consider a rotation on the circle versus a homeomorphism with both periodic and wandering points. Although both have topological entropy equal to zero, the first system is evidently simpler than the second. Labrousse proved that polynomial entropy can distinguish between such systems (see Theorem 1 in~\cite{L}).

Topological and polynomial entropy share several properties: both are conjugacy invariants, they depend only on the topology (not on the specific choice of metric), they satisfy the finite union property, and they both obey a product formula. However, there are important differences as well. Certain properties such as the power formula, the $\sigma$-union property, and the variational principle hold for topological entropy but do not necessarily hold for polynomial entropy (see~\cite{M1,M2}). Moreover, topological entropy is equal to the entropy of the system restricted to the non-wandering set, which is a closed, invariant subset. This is not the case for polynomial entropy. Unlike topological entropy, which captures only the behavior on the non-wandering set, polynomial entropy is sensitive to the wandering part of the system.

In~\cite{DK}, the authors computed the polynomial entropy of the induced map $C(f)$, when $f$ is a circle homeomorphism with finitely many non-wandering points. Those results heavily rely on the coding techniques developed in~\cite{HL} and slightly generalized in~\cite{KP}, which apply only to the case when the non-wandering set is finite (so every non-wandering point is actually periodic).

\begin{thrm}(\cite{DK})\label{thm:DK}
Let $f:[0,1]\to [0,1]$ or $f:\mathbb{S}^1\to \mathbb{S}^1$ be a homeomorphism with a finite non-wandering set. Then $h_{\mathrm{pol}}(C(f))=2$, $h_{\mathrm{pol}}(F_k(f))=k$ and $h_{\mathrm{pol}}(2^{f})=\infty$.\qed
\end{thrm}

In the same paper, we generalised these results for an infinite non-wandering set, in case of $f:[0,1]\to [0,1]$. However, the case of $f:\mathbb{S}^1\to \mathbb{S}^1$ with infinite non-wandering set remained unsolved, for non-trivial cases. In this paper, we compute the polynomial entropy of $C(f)$, for an arbitrary circle homeomorphism $f:\mathbb{S}^1\to \mathbb{S}^1$. Moreover, we generalize the result from~\cite{DK} to the case of the induced map $C_k(f)$.

\begin{thrm}\label{thm:main} Let $f:X\to X$ be a homeomorphism, where $X$ is either $[0,1]$ or $\mathbb{S}^1$ . If $X=\mathbb{S}^1$ and $f$ is conjugate to a rotation, then $h_{\mathrm{pol}}(C_k(f))=0$, otherwise $h_{\mathrm{pol}}(C_k(f))=2k$. If $X=[0,1]$ and $f\neq\mathrm{Id}$ and $f^2\neq\mathrm{Id}$, then $h_{\mathrm{pol}}(C_k(f))=2k$, otherwise $h_{\mathrm{pol}}(C_k(f))=0$.
\end{thrm}

The part concerning the interval is proved in~\cite{KMP}, using a different approach.

\section{Preliminaries}

\subsection{Hyperspaces and induced maps} For a compact metric space $(X,d)$, the hyperspace $2^X$ is the set of all nonempty closed subsets of $X$. The topology on $2^X$ is induced by the Hausdorff metric
$$d_H(A,B):=\inf\{\varepsilon>0\mid A\subset U(B,\varepsilon),\;B\subset U(A,\varepsilon)\},$$
where
$$U(A,\varepsilon):=\{x\in X\mid d(x,A)<\varepsilon\}.$$
The space $2^X$ is compact and the topology induced by $d_H$ is the Vietoris topology.

If $X$ is also connected, one can consider a closed subspace of $2^X$, with the induced metric, namely $C(X)$. The set $C(X)$ of all closed and connected nonempty subsets of $X$ is called the \textit{hyperspace of subcontinua} of $X$ and is also a continuum. For $k\ge 1$, he space $C_k(X)$ is defined as 
$$C_k(X):=\{K\in 2^X\mid K\;\mbox{has at most}\;k\;\mbox{connected components}\}$$ and
$k$-fold symmetric space 
$$F_k(X):=\{K\in2^X\mid |K|\le k\},$$ where $|\cdot|$ denotes the cardinality of a set.

If $f:X\to X$ is continuous, then it induces continuous maps
$$\begin{aligned}&2^f:2^X\to 2^X,\quad &2^f(A):=\{f(x)\mid x\in A\}\\
&C_k(f):C_k(X)\to C_k(X),\quad &C_k(f)(A):=\{f(x)\mid x\in A\}\\
&F_k(f):F_k(X)\to F_k(X),\quad &F_k(f)(A):=\{f(x)\mid x\in A\}.
\end{aligned}$$
If $f$ is a homeomorphism, then so are $2^f$, $C_k(f)$ and $F_k(f)$.

\begin{lem} 
If $f:X\to X$ is an isometry on a compact metric space $(X,d)$, then so is the induced map $C_k(f):C_k(X)\to C_k(X)$.
\end{lem}
\noindent{\it Proof.} Let $A,B\in C_k(X)$ and denote $r=d_H(A,B)$. Since $A\subset U(B,r)$ and $f$ is an isometry, we have that $f(A)\subset U(f(B),r)$. Similarly, we obtain that $f(B)\subset U(f(A),r)$, so $d_H(C_k(f)(A),C_k(f)(B))\leqslant r$.  Let $0<s<r$. We can find $a\in A$ and $b\in B$ such that $d(f(a),f(b))=d(a,b)>s$. Because of this, we cannot have $f(A)\subset U(f(B),s)$, so $d_H(C_k(f)(A),C_k(f)(B))=r=d_H(A,B)$.\qed

\subsection{Polynomial entropy}\label{subsec:entropy}

Suppose that $X$ is a compact metric space, and $f:X\rightarrow X$ is continuous. For each positive integer $n$, denote by $d_n^f(x,y)$ the dynamic metric (induced by $f$ and $d$):
$$
d_n^f(x,y)=\max\limits_{0\leq l\leq n-1}d(f^l(x),f^l(y)).
$$
Fix $Y\subseteq X$. For $\varepsilon>0$, we say that a finite set $E\subset X$ is \textit{$(n,\varepsilon)$-separated} if for every $x,y \in E$ it holds $d_n^f(x,y)\geq \varepsilon$. Let $\mathrm{sep}(n,\varepsilon;Y)$ denote the maximal cardinality of an $(n,\varepsilon)$-separated set $E$, contained in $Y$.

\begin{def}\label{def:pol_ent} The \textit{polynomial entropy} of the map $f$ on the set $Y$ is defined by
$$
h_{\mathrm{pol}}(f;Y)=\lim\limits_{\varepsilon \rightarrow 0}\limsup\limits_{n\rightarrow \infty}\frac{\log \mathrm{sep}(n,\varepsilon;Y)}{\log n}.
$$
\end{def}

We can also define the polynomial entropy as follows. Let $\mathrm{span}(n,\varepsilon;Y)$ denote the minimal number of balls of radius $\varepsilon$ (with respect to $d_n^f$) that cover $Y$ (the centers of the ball are not necessarily elements of $Y$). Denote by $\mathrm{cov}(n,\varepsilon;Y)$ the minimal number of sets $Y_j$ such that the diameters (with respect to $d_n^f$) of $Y_j$ are smaller than $\varepsilon$ and $Y\subseteq\cup_{j=1}^mY_j$.

From the following sequence of inequalities
\begin{equation}\label{eq:span-cov-sep}
\mathrm{span}(n,2\varepsilon;Y)\le \mathrm{cov}(n,\varepsilon;Y)\le \mathrm{sep}(n,\varepsilon;Y)\le \mathrm{cov}(n,\varepsilon/2;Y)\le \mathrm{span}(n,\varepsilon/2;Y)
\end{equation} 
we conclude that
$$h_{\mathrm{pol}}(f;Y)=\lim\limits_{\varepsilon \rightarrow 0}\limsup\limits_{n\rightarrow \infty}\frac{\log \mathrm{span}(n,\varepsilon;Y)}{\log n}=\lim\limits_{\varepsilon \rightarrow 0}\limsup\limits_{n\rightarrow \infty}\frac{\log \mathrm{cov}(n,\varepsilon;Y)}{\log n}.$$

If $X=Y$ we abbreviate $h_{\mathrm{pol}}(f):=h_{\mathrm{pol}}(f;X)$. We list some properties of the polynomial entropy that are important for our computations (for proofs see Propositions $1-4$ in \cite{M2}):
\begin{itemize}
\item[(1)] $h_{\mathrm{pol}}(f^k)=h_{\mathrm{pol}}(f)$, for any $k\geqslant 1$.
\item[(2)] If $Y\subset X$ is a closed, $f$-invariant set, then $h_{\mathrm{pol}}(f;Y)=h_{\mathrm{pol}}(f|_Y)$.
\item[(3)] If $Y=\bigcup_{j=1}^mY_j$ where each $Y_j$ are $f$-invariant, then $h_{\mathrm{pol}}(f;Y)=\max\{h_{\mathrm{pol}}(f;Y_j)\mid j=1,\ldots,m\}$\label{finite-union}.\label{(3)}
\item[(4)] If $f:X\to X$, $g:Y\to Y$ and $f\times g:X\times Y\to X\times Y$ is defined as $f\times g (x,y):=(f(x),g(y))$, then $h_{\mathrm{pol}}(f\times g)=h_{\mathrm{pol}}(f)+h_{\mathrm{pol}}(g)$.
\item[(5)] $h_{\mathrm{pol}}(f)$ does not depend on a metric but only on the induced topology.
\item[(6)] $h_{\mathrm{pol}}(\cdot)$ is a \textit{conjugacy invariant} (meaning if $f:X\to X$, $g:X'\to X'$, $\varphi:X\to X'$ is a homeomorphism of compact spaces and $g\circ\varphi=\varphi\circ f$, then $h_{\mathrm{pol}}(f)=h_{\mathrm{pol}}(g)$).
\item[(7)]\label{page:properties} If $f:X\to X$ and $g:X'\to X'$ are \textit{semi-conjugated}, meaning that $\varphi:X\to X'$ is a continuous surjective map of compact spaces and $g\circ\varphi=\varphi\circ f$, then $h_{\mathrm{pol}}(f)\geqslant h_{\mathrm{pol}}(g)$.
\end{itemize}

A set $A\subset X$ is \textit{wandering} if $f^n(A)\cap A=\emptyset$, for all $n\geqslant 1$. A point $p\in X$ is \textit{wandering} if there exists a wandering neighbourhood $U\ni p$. \label{wan-pt}

A point that is not wandering is said to be \textit{non-wandering}. We denote the set of all non-wandering points by $\mathrm{NW}(f)$. The set $\mathrm{NW}(f)$ is closed and $f$-invariant. Also, we denote the set of all fixed points by $\mathrm{Fix}(f)$ and the set of all periodic points by $\mathrm{Per}(f)$.

\subsection{Circle homeomorphisms}

Starting from the dynamical characteristics of rational and irrational rotation on the circle $\mathbb{S}^1$, one can study the behaviour of arbitrary orientation-preserving homeomorphisms on $\mathbb{S}^1$ by means of \textit{Poincar\'e rotation number}. As it is expected, homeomorphisms with a rational rotation number have periodic points and simpler dynamics in general, while that is not the case with homeomorphisms with an irrational rotation number.\\

\noindent Let $f:\mathbb{S}^1\to\mathbb{S}^1$ be an orientation-preserving homeomorphism. Recall the classification via Poincar\'e rotation number (denoted by $\rho(f)$):
\begin{itemize}
\item If $\rho(f)$ is rational, $\rho(f)=p/q$, then $\mathrm{NW}(f)=\mathrm{Per}(f)$. Also, there exists a periodic point and all the periodic points have the same period - $q$. If there exist no non-periodic points, then all the points are periodic, with the period $q$, so $f^q=\mathrm{Id}$. The latter case occurs if and only if $f$ is conjugated to the rotation by $\rho(f)$.\label{page:rot-no}

\item If $\rho(f)\notin\mathbb{Q}$, then there are two possibilities. If $f$ is topologically transitive, then $f$ is conjugate to the rotation by $\rho(f)$. If not, the set $\mathrm{NW}(f)$ is homeomorphic to a Cantor set. In both cases, $\omega$ and $\alpha$-limit sets do not depend on a particular point, and coincide with the non-wandering set.
\end{itemize}

For more details on the rotation number of circle homeomorphisms and proofs of the stated properties, see \cite{T}.

\section{Main results and proofs}

We present the proof of the main Theorem \ref{thm:main} in several steps.

We first prove the following general lemma.
\begin{lem}\label{lem:gen} For any continuous self map $f$ of a compact metric space $X$
$$\mathrm{span}(n,\varepsilon,C_k(f))\le k\cdot\mathrm{span}(n,\varepsilon,C(f))^k,$$ and therefore
$$h_{\mathrm{pol}}(C_k(f))\le k\cdot h_{\mathrm{pol}}(C(f)) .$$
\end{lem}
\noindent{\it Proof.} Let $A_1,\ldots,A_N$ be points in $C(X)$ such that the balls of diameter $\varepsilon$ with respect to the metric $d_n^{\,C(f)}$, centered at $A_j$, cover $C(X)$ and let $N$ be minimal such that this is true. For $l\in\{1,\ldots,k\}$ and $j_i\in\{1,\ldots,N\}$. Define
$$D_{j_1\ldots j_l}:=A_{j_1}\cup\ldots\cup A_{j_l}.$$ It is obvious that $D_{j_1\ldots j_l}\in C_k(X)$. The balls of diameter $\varepsilon$ with respect to the metric $d_n^{\,C_k(f)}$, centered at $D_{j_1\ldots j_l}$, cover $C_k(X)$. Indeed, for any $K\in C_k(X)$, $K=K_1\cup\ldots\cup K_m$, where $K_j\in C(X)$ and $m\le k$. There exists $p_j$ such that
$$d_n^{\,C(f)}(A_{p_j},K_j)<\varepsilon$$ for every $j\in\{1,\ldots,m\}$. We abbreviate the notation $(d_H)_n^{\,C(f)}$ to $d_n^{\,C(f)}$ and we will use this abbreviation throughout the rest of the paper. Therefore
$$
d_n^{\,C_k(f)}(D_{p_1\ldots p_m},K)<\varepsilon.$$ Now we compute the cardinality of the set 
$$\mathcal{D}:=\left\{D_{j_1\ldots j_l}\mid l\in\{1,\ldots,k\},\,  j_i\in\{1,\ldots,N\}\right\}.$$ We have
$$|\mathcal{D}|\le N+{N\choose 2}+\ldots {N\choose k}\le N+N^2+\ldots+N^k\le k\cdot{N^k}.$$\qed

\begin{prop}\label{prop:le}
For any continuous $f:\mathbb{S}^1\to\mathbb{S}^1$ we have
$$h_{\mathrm{pol}}(C_k(f))\leqslant 2k\cdot h_{\mathrm{pol}}(f).$$
\end{prop}
\noindent{\it Proof.} For the first part of the proof, we follow the proof of Lemma 4.1 in \cite{HCS}, which is a proposition regarding the topological entropy. Let $\varepsilon>0$ and $n\in\mathbb{N}$. Let $N=\mathrm{cov}\left(n,\dfrac{\varepsilon}{2},\mathbb{S}^1\right)$ and suppose that $I_0,I_1,\ldots,I_{N-1}$ are arcs (numbered in counterclockwise direction) with $d_n^f$-diameter less than $\varepsilon/2$ that cover $\mathbb{S}^1$. We will now construct a cover of $C(\mathbb{S}^1)$ with no more than $N^2$ elements whose $d_{n}^{\,C(f)}$-diameter is less than $\varepsilon$. For every $I_j$, we define families of arcs $M_k(I_j)$, for all $k=0,\ldots,N-1$ in the following way (we use the convention $I_j=I_{j \,\mathrm{mod} \ N}$):

$$\begin{aligned}
&M_0(I_j):=C(I_j)\\
&M_1(I_j):=\{I\in C(I_j\cup I_{j-1})\ \big{|}\ I\supset I_j\cap I_{j-1}\}\\
&M_2(I_j):=\{I\in C(I_{j+1}\cup I_j\cup I_{j-1})\ \big{|}\ I\supset I_j\}\\
&M_3(I_j):=\{I\in C(I_{j+1}\cup I_j\cup I_{j-1}\cup I_{j-2})\ \big{|}\ I\supset I_j\cup I_{j-1}\}\\
&M_4(I_j):=\{I\in C(I_{j+2}\cup I_{j+1}\cup I_j\cup I_{j-1}\cup I_{j-2})\ \big{|}\ I\supset I_{j+1}\cup I_j\cup I_{j-1}\}\\
&\ldots\ldots\ldots\\
&M_{N-1}(I_j):=\{I\in C(\mathbb{S}^1)\ \big{|}\ I\supset \mathbb{S}^1\setminus I_{j+N/2}\}.
\end{aligned}$$

It is clear that $d^{\,C(f)}_n$-diameter of each $M_k(I_j)$ is less than $\varepsilon$. Also, the collection $\{M_k(I_j)\ \big{|} \ k=0,\ldots,N-1, j=0,\ldots,N-1\}$ forms a covering of $C(\mathbb{S}^1)$. Indeed, take any $[a,b]\in C(\mathbb{S}^1)$ and let $j_1$ be the least index with $a\in I_{j_1}$. If $b\in I_{j_1}$, we have $[a,b]\in M_0$. Otherwise, let $j_2>j_1$ be the least index with $b\in I_{j_2}$. We see that $[a,b]\in M_k$, where $k:=j_2-j_1$. 

We now have that $\mathrm{cov}\left(n,\varepsilon,C(f)\right)\leqslant N^2$. From Lemma~\ref{lem:gen} and the sequence of inequalities~(\ref{eq:span-cov-sep}), we conclude
$$\mathrm{span}(n,2\varepsilon,C_k(f))\le k\cdot\mathrm{span}(n,2\varepsilon,C(f))^k\le k\cdot \mathrm{cov}(n,\varepsilon,C(f))^k\le k\cdot N^{2k},$$
and consequently that $$h_{\mathrm{pol}}(C_k(f))\leqslant 2k\cdot h_{\mathrm{pol}}(f).$$

\qed

Using that $h_{\mathrm{pol}}(f)\leqslant1$ for any  orientation-preserving homeomorphism $f:\mathbb{S}^1\to\mathbb{S}^1$ (see \cite{L}), we can deduce the following corollary.
\begin{cor}\label{cor:ge}
For any orientation-preserving homeomorphism $f:\mathbb{S}^1\to\mathbb{S}^1$ we have
$$h_{\mathrm{pol}}(C_k(f))\leqslant 2k.$$
\end{cor}

\begin{prop}\label{prop:geQ} Let $f:\mathbb{S}^1\to\mathbb{S}^1$ be an orientation-preserving homeomorphism that possesses a wandering point. Then $h_{\mathrm{pol}}(C_k(f))\geqslant 2k$.
\end{prop}
\noindent{\it Proof.} In~\cite{DKL} we proved that, under the assumption of Proposition~\ref{prop:geQ}, $h_{\mathrm{pol}}(F_m(f))\ge m$. Since $C_k(f)$ is a homeomorphism, the sets $\{\mathbb{S}^1\}$ and $Y:=C_k(\mathbb{S}^1)\setminus\{\mathbb{S}^1\}$ are $C(f)$-invariant and their union is the whole space $C_k(\mathbb{S}^1)$. Therefore 
$$\begin{aligned}
&h_{\mathrm{pol}}(C_k(f))=\max\{h_{\mathrm{pol}}(C_k(f)|_Y),h_{\mathrm{pol}}\left(C_k(f)|_{\{\mathbb{S}^1\}}\right)\}=\\
&\max\{h_{\mathrm{pol}}(C_k(f)|_Y),0\}=h_{\mathrm{pol}}(C_k(f)|_Y).\end{aligned}$$
Consider the following continuous surjection:
$$\pi:Y\to F_{2k}(\mathbb{S}^1),\quad\pi:K\mapsto\partial(K)$$
where $\partial(K)$ denotes the boundary of $K$. It is obvious that 
$$\pi\circ C_k(f)|_Y=F_{2k}(f)\circ\pi,$$
so 
\begin{equation}\label{eq:ineq}
h_{\mathrm{pol}}(C_k(f))=h_{\mathrm{pol}}(C_k(f)|_Y)\ge h_{\mathrm{pol}}(F_{2k}(f))\ge 2k.
\end{equation}
Note that the latter inequality does not immediately follow from property (7) of the polynomial entropy stated on page~\pageref{page:properties}, since the space $Y$ is not compact. However, $\pi$ is a Lipschitz map with Lipschitz constant $1$. Thus, if the set $E=\{A_1,\ldots,A_N\}\subset F_{2k}(\mathbb{S}^1)$ is an $(n,\varepsilon)$-separated set such that $\mathrm{sep}(n,\varepsilon,F_{2k}(f))=|E|$, let $K_j\in C(\mathbb{S}^1)$ be such that $\pi(K_j)=A_j$. We have
$$\begin{aligned}
&d(C_k(f)^l(K_i),C_k(f)^l(K_j))\ge d(\pi(C_k(f)^l(K_i)),\pi(C_k(f)^l(K_j)))=\\
&d(F_{2k}(f)^l(\pi(K_i)),F_{2k}(f)^l(\pi(K_j)))=d(F_{2k}(f)^l(A_i),F_{2k}(f)^l(A_j)),
\end{aligned}$$ so
$$d_n^{\,C_k(f)}(K_i,K_j)\ge d_n^{F_{2k}(f)}(A_i,A_j),$$ and therefore the set 
$$\widetilde{E}:=\{K_1,\ldots,K_N\}\subset C(\mathbb{S}^1)$$ is also an $(n,\varepsilon)$-separated set. We conclude 
$$ \mathrm{sep}(n,\varepsilon,C_k(f))\ge \mathrm{sep}(n,\varepsilon,F_{2k}(f)),$$ so~(\ref{eq:ineq}) holds.

\qed

\begin{prop}\label{prop:geQ-I} Let $f:[0,1]\to[0,1]$ be a diffeomorphism such that $f\neq\mathrm{Id}$ and $f^2\neq\mathrm{Id}$. Then $h_{\mathrm{pol}}(C_k(f))\geqslant 2k$.
\end{prop} 
\noindent{\it Proof.} We can assume that $f$ is increasing, since, otherwise $f^2$ is. If $f\neq\mathrm{Id}$ there exists $x\in [0,1]$ such that $f(x)>x$ or $f(x)<x$. In both cases $x$ is wandering, so by a result from~\cite{DKL}, $h_{\mathrm{pol}}(F_m(f))\ge m$. The rest of the proof is the same as the proof of Proposition~\ref{prop:geQ}, except that here we take $Y$ to be the whole $C_k([0,1])$ and the boundary of $K$ is taken in $\mathbb{R}$.\qed

\vspace{7mm}

We can now proceed to prove the main theorem \ref{thm:main}.\\

\noindent{\it Proof of Theorem~\ref{thm:main}.} First we deal with the circle. We can assume that $f$ is orientation-preserving; otherwise, $f^2$ is orientation-preserving, and 
$$h_{\mathrm{pol}}\left(C_k\left(f^2\right)\right)=h_{\mathrm{pol}}\left(C_k(f)^2\right)=h_{\mathrm{pol}}(C_k(f)).$$ 

If $f$ is conjugate to a rotation $\rho_\alpha$, then $C_k(f)$ is conjugate to $C_k(\rho_\alpha)$, which is an isometry (since $\rho_\alpha$ is), so $h_{\mathrm{pol}}(C_k(f))=0$. 

By Proposition~\ref{prop:le} we have $h_{\mathrm{pol}}(C_k(f))\leqslant 2k$. 

If the rotation number of $f$ is rational, then, by the discussion on page~\pageref{page:rot-no} we conclude that $f$ has a wandering point. Hence, by Proposition~\ref{prop:geQ}, we obtain $h_{\mathrm{pol}}(C_k(f))\geqslant 2k$. 

If the rotation number of $f$ is irrational, we consider two cases: the non-wandering set of $f$ is either $\mathbb{S}^1$ or homeomorphic to the Cantor set. In the first case, $f$ is conjugate to a rotation, so we are done. In the second case, we again deduce that $f$ has a wandering point, and therefore, by Proposition~\ref{prop:geQ}, $h_{\mathrm{pol}}(C_k(f))\geqslant 2k$.

The case of the interval follows imediately from Theorem~\ref{thm:DK}, Lemma~\ref{lem:gen} and Proposition~\ref{prop:geQ-I}.
\qed

\smallskip

\vspace{5mm}

\textbf{Acknowledgements}: This work is partially supported by the Ministry of Education, Science and Technological Developments of Republic of Serbia: grant number 451-03-47/2023-01/200104 with Faculty of Mathematics, University of Belgrade. \\

\end{document}